\documentclass[12pt,leqno]{article}
\usepackage{amsfonts}
\usepackage{amsmath, amsthm, amsfonts, amssymb, color}
\usepackage{mathrsfs}
\usepackage{color}
\newtheorem{thm}{Theorem}[section]

\newtheorem{rem}[thm]{Remark}
\theoremstyle{definition}

\newcommand{\scr}[1]{\mathscr #1}
\definecolor{wco}{rgb}{0.5,0.2,0.3}

\numberwithin{equation}{section} \theoremstyle{remark}

\newcommand{\ua}{\uparrow}

\title{{\bf Comparison Theorem for Path Dependent SDEs Driven by $G$-Brownian Motion}\footnote{Supported in
 part by  NNSFC (11801406).} }
\author{
{\bf     Xing Huang $^{1,a)}$, Fen-Fen Yang $^{1,b)}$   }\\
\footnotesize{  1)Center for Applied Mathematics, Tianjin University, Tianjin 300072, China}\\
\footnotesize{ $^{a)}$xinghuang@tju.edu.cn, $^{b)}$yangfenfen@tju.edu.cn}
}
\begin{document}
\allowdisplaybreaks
\def\R{\mathbb R}  \def\ff{\frac} \def\ss{\sqrt} \def\B{\mathbf
B} \def\W{\mathbb W}
\def\N{\mathbb N} \def\kk{\kappa} \def\m{{\bf m}}
\def\ee{\varepsilon}\def\ddd{D^*}
\def\dd{\delta} \def\DD{\Delta} \def\vv{\varepsilon} \def\rr{\rho}
\def\<{\langle} \def\>{\rangle} \def\GG{\Gamma} \def\gg{\gamma}
  \def\nn{\nabla} \def\pp{\partial} \def\E{\mathbb E}
\def\d{\text{\rm{d}}} \def\bb{\beta} \def\aa{\alpha} \def\D{\scr D}
  \def\si{\sigma} \def\ess{\text{\rm{ess}}}
\def\beg{\begin} \def\beq{\begin{equation}}  \def\F{\scr F}
\def\Ric{\text{\rm{Ric}}} \def\Hess{\text{\rm{Hess}}}
\def\e{\text{\rm{e}}} \def\ua{\underline a} \def\OO{\Omega}  \def\oo{\omega}
 \def\tt{\tilde} \def\Ric{\text{\rm{Ric}}}
\def\cut{\text{\rm{cut}}} \def\P{\mathbb P} \def\ifn{I_n(f^{\bigotimes n})}
\def\C{\scr C}      \def\aaa{\mathbf{r}}     \def\r{r}
\def\gap{\text{\rm{gap}}} \def\prr{\pi_{{\bf m},\varrho}}  \def\r{\mathbf r}
\def\Z{\mathbb Z} \def\vrr{\varrho} \def\ll{\lambda}
\def\L{\scr L}\def\Tt{\tt} \def\TT{\tt}\def\II{\mathbb I}
\def\i{{\rm in}}\def\Sect{{\rm Sect}}  \def\H{\mathbb H}
\def\M{\scr M}\def\Q{\mathbb Q} \def\texto{\text{o}} \def\LL{\Lambda}
\def\Rank{{\rm Rank}} \def\B{\scr B} \def\i{{\rm i}} \def\HR{\hat{\R}^d}
\def\to{\rightarrow}\def\l{\ell}\def\iint{\int}
\def\EE{\scr E}\def\Cut{{\rm Cut}}
\def\A{\scr A} \def\Lip{{\rm Lip}}
\def\BB{\scr B}\def\Ent{{\rm Ent}}\def\L{\scr L}
\def\R{\mathbb R}  \def\ff{\frac} \def\ss{\sqrt} \def\B{\mathbf
B}
\def\N{\mathbb N} \def\kk{\kappa} \def\m{{\bf m}}
\def\dd{\delta} \def\DD{\Delta} \def\vv{\varepsilon} \def\rr{\rho}
\def\<{\langle} \def\>{\rangle} \def\GG{\Gamma} \def\gg{\gamma}
  \def\nn{\nabla} \def\pp{\partial} \def\E{\mathbb E}
\def\d{\text{\rm{d}}} \def\bb{\beta} \def\aa{\alpha} \def\D{\scr D}
  \def\si{\sigma} \def\ess{\text{\rm{ess}}}
\def\beg{\begin} \def\beq{\begin{equation}}  \def\F{\scr F}
\def\Ric{\text{\rm{Ric}}} \def\Hess{\text{\rm{Hess}}}
\def\e{\text{\rm{e}}} \def\ua{\underline a} \def\OO{\Omega}  \def\oo{\omega}
 \def\tt{\tilde} \def\Ric{\text{\rm{Ric}}}
\def\cut{\text{\rm{cut}}} \def\P{\mathbb P} \def\ifn{I_n(f^{\bigotimes n})}
\def\C{\scr C}      \def\aaa{\mathbf{r}}     \def\r{r}
\def\gap{\text{\rm{gap}}} \def\prr{\pi_{{\bf m},\varrho}}  \def\r{\mathbf r}
\def\Z{\mathbb Z} \def\vrr{\varrho} \def\ll{\lambda}
\def\L{\scr L}\def\Tt{\tt} \def\TT{\tt}\def\II{\mathbb I}
\def\i{{\rm in}}\def\Sect{{\rm Sect}}  \def\H{\mathbb H}
\def\M{\scr M}\def\Q{\mathbb Q} \def\texto{\text{o}} \def\LL{\Lambda}
\def\Rank{{\rm Rank}} \def\B{\scr B} \def\i{{\rm i}} \def\HR{\hat{\R}^d}
\def\to{\rightarrow}\def\l{\ell}
\def\8{\infty}\def\I{1}\def\U{\scr U}
\maketitle

\begin{abstract} Sufficient and necessary conditions  are presented for the comparison theorem of  path dependent $G$-SDEs. Different from  the corresponding study in path independent $G$-SDEs, a probability method is applied to prove these results. Moreover, the results extend the ones in the linear expectation case.
\end{abstract} \noindent
 \vskip 2cm

\section{Introduction}

The order preservation of stochastic processes is an important property for one to compare a complicated process with simpler ones, and a result to ensure this property is called $``$comparison theorem" in the literature. There are two different type order preservations, one is in the distribution (weak) sense and the other is in the pathwise (strong) sense, where the latter implies the former.

In the linear expectation frame, the weak order preservation has been investigated in \cite{CW,W1,W} and references within. There are also lots of results on the strong order preservation, see, for instance,
   \cite{BY, GD, HW, IW, M, O, PY, PZ, YMY, Z} and references therein for comparison theorems on forward/backward SDEs (stochastic differential equations), with jumps and/or with  memory. Recently, the first author and his co-authors extend the results in \cite{HW} to the path-distribution dependent case, one can refer to \cite{HLW} for more details.

On the other hand, there are some results on the comparison theorem for $G$-SDEs, see \cite{Lin,luo,LW}. Some sufficient condition is presented in \cite[Theorem 7.1]{Lin} for comparison theorem of one-dimensional $G$-SDEs. In \cite{LW}, the authors obtain the sufficient and necessary conditions for comparison theorem by the viability property of SDEs, which is equivalent to the fact that the square of the distance to the constraint set is a viscosity supersolution to the associated
 Hamilton-Jacobi-Bellman equation, see \cite[Theorem 2.5]{LW} and references therein for more details.

The aim of this paper is to present sufficient and necessary conditions of the order preservations for path dependent $G$-SDEs and we provide a probability method to prove them. The result extends the ones in \cite{HW} when the noise is standard Brownian motion. We will adopt the method in \cite{HW} to  complete the proof. However, some essential work needs to been done since the quadratic variation process $\<B\>$ of the $G$-Brownian motion $B$ is not determined under $G$-expectation. More precisely, we need to treat $\int_0^\cdot\<h(s), \d \<B\>(s)\>-2\int_0^\cdot G(h(s))\d s$ which is well known as a non-increasing $G$-martingale. This is quite different from the linear expectation case. Moreover, in the proof of necessary condition of the comparison theorem, we will use the representation theorem \eqref{rep2} below of the $G$-expectation introduced in \cite{15,HP,song}, by which the order preservation under $G$-expectation implies that in linear expectation case. Then the existed result in \cite{HW} can be applied to prove the necessary condition on diffusion coefficients.

Before moving on, we recall some basic facts on $G$-expectation and $G$-Brownian motion in the following section.

\section{$G$-Expectation and $G$-Brownian motion}
 Let $\Omega=C_0([0,\infty);\mathbb{R}^m)$, the $\mathbb{R}^m$-valued and continuous functions on $[0,\infty)$ vanishing at zero, equipped with the metric
$$\rho(\omega^1,\omega^2)=\sum_{n=1}^\infty\frac{1}{2^n}\left[\max_{t\in[0,n]}|\omega^1_t-\omega^2_t|\wedge1 \right],\ \ \omega^1,\omega^2\in\Omega.$$
%Let $T>0$ and $\Omega_T=C_0([0,T];\mathbb{R}^d)$, the $\mathbb{R}^d$-valued and continuous functions on $[0,T]$ vanishing at zero.
For any $T>0$, set
 $$L_{ip}(\Omega_T) =\{\omega\to \varphi(\omega_{t_1}, \cdot \cdot \cdot, \omega_{t_n}):n\in \mathbb{N}^+, t_1,\cdot \cdot \cdot, t_n\in [0,T],\varphi \in C_{b,lip}((\mathbb{R}^{m})^n)\},$$
 and $$L_{ip}(\Omega)=\bigcup_{T>0}L_{ip}(\Omega_T),$$
where  $C_{b,lip}((\mathbb{R}^{m})^n)$ denotes  the set of bounded and Lipschitz continuous functions on $((\mathbb{R}^{m})^n)$. We denote by $|A|^2=\|A\|_{HS}^2$ for any matrix $A$. For two $m\times m$ matrices $M$ and $\bar{M}$, define $$\<M,\bar{M}\>=\sum_{k,l=1}^mM _{kl}\bar{M} _{kl}.$$
Let $\mathbb{M}^m$ be the collection of all $m\times m$ matrices and $\mathbb{S}^m$ ($\mathbb{S}_+^m$) be the set of the symmetric (symmetric and positive definite) ones in $\mathbb{M}^m$.
%Fix $\underline{\sigma}, \bar{\sigma}\in \mathbb{S}_+^d$ with $\underline{\sigma}<\bar{\sigma}$ and
Fix two positive constants $\underline{\sigma}<\bar{\sigma}$ and  define
\begin{equation}\label{G(A)}
 G(A):=\frac{1}{2}\sup _{\gamma\in \mathbb{S}_+^m\bigcap[\underline{\sigma}^2\textbf{I}_{m\times m}, \bar{\sigma}^2\textbf{I}_{m\times m}]}\<\gamma,A\>, \ A\in\mathbb{S}^m.
\end{equation}
It is  not difficult to see that $G$ has the following properties:
\begin{enumerate}
\item[(a)] (Positive homogeneity) $G(\lambda A)=\lambda G(A)$, $\lambda\geq 0,  A\in\mathbb{S}^m$.
\item[(b)] (Sub-additivity) $G(A+\bar{A})\leq G(A)+G(\bar{A})$, $G(A)-G(\bar{A})\leq G(A-\bar{A})$,\ \   $A,\bar{A}\in\mathbb{S}^m$.
\item[(c)] $|G(A)|\leq \frac{1}{2}|A|\sup _{\gamma\in \mathbb{S}_+^m\bigcap[\underline{\sigma}^2\textbf{I}_{m\times m}, \bar{\sigma}^2\textbf{I}_{m\times m}]}|\gamma|=\frac{1}{2}|A|\sqrt{m}\bar{\sigma}^2.$
\item[(d)] $G(A)-G(\bar{A})\geq \frac{\underline{\sigma}^2}{2} \mbox{trace}[A-\bar{A}],\ A\geq \bar{A}, A,\bar{A} \in \mathbb{S}^m.$
%where $\lambda_0(\underline{\sigma}^2)>0$ is the minimal eigenvalue of $\underline{\sigma}^2$.
\end{enumerate}
\begin{rem}\label{coG}
(b) and (c) imply that $G$ is continuous.
\end{rem}
 Let $\bar{\E}^G$ be the nonlinear expectation on $\Omega$ such that coordinate process $(B(t))_{t\geq 0}$, i.e. $B(t)(\omega)=\omega_t, \omega\in \Omega$, is an $m$-dimensional $G$-Brownian motion on $(\Omega, L_G^1(\Omega),\bar{\E}^G)$, where $L_G^1(\Omega)$ is the completion of $L_{ip}(\Omega)$  under the norm $\bar{\E}^G|\cdot|$. One can refer to \cite{song} for details on the construction of $\bar{\E}^G$. For any $p\geq 1$, let $L_G^p(\Omega)$ be the completion of $L_{ip}(\Omega)$  under the norm $(\bar{\E}^G|\cdot|^p)^{\frac{1}{p}}$. Similarly, we can define $L_G^p(\Omega_T)$ for any $T>0$.

Let
\begin{align*}
M_{G}^{p,0}([0,T])&=
\Big\{\eta_t:=\sum_{j=0}^{N-1} \xi_{j} I_{[t_j, t_{j+1})}(t);
~\xi_{j}\in L_{G}^p(\Omega_{t_{j}}),
 N\in\mathbb{N}^+,\\
 & ~~~\qquad \qquad0=t_0<t_1<\cdots <t_N=T \Big\},
\end{align*}
and $M_G^p([0,T])$ be the completion of $M_G^{p,0}([0,T])$ under the norm
$$\|\eta\|_{M_G^p([0,T])}:=\left(\bar{\E}^G\int_{0}^{T}|\eta_{t}|^p\d t\right)^{\frac{1}{p}}.$$
%Moreover, let $$ M_G^2([0,T])^d=\left\{X=(X^1,X^2,\cdots,X^d), X^i\in M_G^2([0,T]),1\leq i\leq d\right\}.$$
Let $\mathcal{M}$ be the collection of all probability measures on  $(\Omega, \B(\Omega))$.
According to \cite{15,HP},    there exists a weakly compact subset $\mathcal{P}\subset \mathcal{M}$   such that
%\begin{align}\label{rep1}\bar{\E}^G[X]=\sup_{P\in \mathcal{P}}\mathbb{E}_P[X], \ X\in L_{ip}(\Omega),\end{align}
%and then
\begin{align}\label{rep}\bar{\E}^G[X]=\sup_{P\in \mathcal{P}}\mathbb{E}_P[X], \ X\in L_G^1(\Omega),\end{align}
where $\mathbb{E}_P$ is the linear expectation under probability measure $P \in \mathcal{P}$. $\mathcal{P}$ is called a set that represents $\bar{\E}^G$.
In fact, let $W^0$ be an $m$-dimensional Brownian motion on a complete filtered probability space $(\hat{\OO},\{\F_t\}_{t\geq 0},\P)$, and define 
\begin{align*}\mathbb{H}:=\{\theta: \ \ &\theta\ \ \text{is an}\ \ \mathbb{M}^m\text{-valued progressively measurable}\\ 
&\text{stochastic process},\ \ 
\theta_s\theta_s^\ast\in[\underline{\sigma}^2\textbf{I}_{m\times m}, \bar{\sigma}^2\textbf{I}_{m\times m}], \ \ s\geq 0\},
 \end{align*}
 here $(\cdot)^\ast$ stands for the transpose of a matrix.
 For any $\theta\in \mathbb{H}$, define $\P_{\theta}$ as the law of $\int_0^\cdot \theta_s \d W^0_s$. Then by \cite{15,HP}, we can take $\mathcal{P}= \{\P_\theta, \theta\in\mathbb{H}\}$, i.e.
\begin{align}\label{rep2}
\bar{\E}^G[X]=\sup_{\theta\in \mathbb{H}}\mathbb{E}_{\P_\theta}[X], \ X\in L_G^1(\Omega).
\end{align}
The associated  Choquet capacity to $\bar{\E}^G$ is defined by
$$\mathcal{C}(A)=\sup_{P\in \mathcal{P}}P(A), \ A\in \B(\Omega).$$
A set $A\in\B(\Omega)$ is called polar if $\mathcal{C}(A)=0$,  and we say that a property   holds  $\mathcal{C}$-quasi-surely ($\mathcal{C}$-q.s.)
if it holds outside a polar set, see \cite{15} for more details on capacity.

Finally, letting $\<B\>$ be the quadratic variation process of $B$, then by property (d) and \cite[Chapter III, Corollary 5.7]{peng4}, we have $\mathcal{C}$-q.s.
\begin{align}\label{B}
\underline{\sigma}^2\textbf{I}_{m\times m}<\frac{\d}{\d t}\langle B\rangle(t)\leq\bar{\sigma}^2\textbf{I}_{m\times m}.
\end{align}

\section{Main Results}
Let $r_0\geq 0$ be a constant and   $d\ge 1$ be a natural number.
$\C=C([-r_0,0];\R^d)$ is equipped with uniform norm $\|\cdot\|_\infty$.
For any continuous map $f: [-r_0,\infty)\to \R^d$  and
$t\ge 0$,  let  $f_t\in\C$ be such that $f_t(s)=f(s+t)$ for $s\in
[-r_0,0]$. We call $(f_t)_{t\ge 0}$ the segment of $(f(t))_{t\ge -r_0}.$

Consider the following path dependent SDEs:
\beq\label{E1} \beg{cases} \d X(t)= b(t,X_t)\,\d t+ \<h(t,X_t),\d \<B\>(t)\> +\sigma(t,X_t)\,\d B(t),\\
\d \bar{X}(t)= \bar{b}(t,\bar{X}_t)\,\d t+ \<\bar{h}(t,\bar{X}_t),\d \<B\>(t)\>+\bar{\sigma}(t,\bar{X}_t)\,\d B(t),\end{cases}\end{equation}
where
\begin{align*}&b,\bar{b}: [0,\infty)\times \C\to \R^d;\ \ h,\bar{h}: [0,\infty)\times \C\to (\R^m\otimes\R^m)^d;\\
&\si,\bar{\sigma}: [0,\infty)\times \C\to \R^d\otimes\R^m
\end{align*}
are   measurable.

Without loss of generality, we assume that for any $i=1,\cdots ,d$, $h^i$ and $\bar{h}^i$ are symmetric. Otherwise, we can replace $h^i$ and $\bar{h}^i$ by $\frac{h^i+(h^i)^\ast}{2}$ and $\frac{\bar{h}^i+(\bar{h}^i)^\ast}{2}$ respectively to symmetrize them.

%For any $s\ge0$ and $\F_s$-measurable $\C$-valued random variables $\xi,\bar\xi$,
For any $s\ge0$ and $\xi,\bar\xi\in\C$,
a solution to (\ref{E1}) for $t\ge s$ with $(X_s,\bar X_s)= (\xi,\bar\xi)$ is a  continuous process $(X(t),\bar X(t))_{t\ge s}$  such that for all $t\ge s,$
\beg{equation*}\beg{split} &X(t) = \xi(0)+ \int_s^t  b(r,X_r)\d r+\int_s^t \<h(r,X_r),\d \<B\>(r)\> +\int_s^t \si(r, X_r)\d B(r ),\\
&\bar X(t) = \bar\xi(0)+ \int_s^t \bar b(r,\bar X_r)\d r+\int_s^t \<\bar{h}(r,\bar{X}_r),\d \<B\>(r)\>+ \int_s^t\bar \si(r,\bar X_r)\d B(r) ,\end{split}\end{equation*} where  $(X_t, \bar X_t)_{t\ge s}$ is the segment process of $(X(t), \bar X(t))_{t\ge s-r_0}$ with
$(X_s,\bar X_s)= (\xi,\bar\xi)$.
Throughout the paper, we make   the following assumptions.

\beg{enumerate} \item[{\bf (H1)}]
There exists an increasing function $\aa: \R_+\to \R_+$     such that for any $ t\ge 0, \xi,\eta\in \C$,
 \begin{align*}
&|b(t,\xi)- b(t,\eta)|^2+|\bar{b}(t,\xi)- \bar{b}(t,\eta)|^2+|h(t,\xi)- h(t,\eta)|^2+|\bar{h}(t,\xi)- \bar{h}(t,\eta)|^2\\
&+\|\si(t,\xi)- \si(t,\eta)\|_{HS}^2+\|\bar{\si}(t,\xi)- \bar{\si}(t,\eta)\|_{HS}^2\le \aa(t)\|\xi-\eta\|_{\infty}^2.
\end{align*}
\item[{\bf (H2)}] There exists an increasing function $K: \R_+\to\R_+   $ such that
 \begin{align*}
  &|b(t,0)|^2+ |\bar b(t,0)|^2+|h(t,0)|^2+ |\bar h(t,0)|^2+\|\si(t,0)\|_{HS}^{2}+\|\bar\si(t,0)\|^{2}_{HS}
\le K (t),\ \ t\ge 0.
 \end{align*}
\end{enumerate}
\begin{rem}
According to \cite[Lemma 2.1]{HH}, under {\bf(H1)}-{\bf(H2)}, for any $s\ge 0$ and $\xi, \bar\xi\in \C$,
%$\xi, \bar\xi\in L^2_G(\OO_s\to\C)$,
 the equation  (\ref{E1}) has a unique solution denoted by $\{X(s,\xi;t), \bar X(s,\bar \xi;t)\}_{t\ge s-r_0}$    with $X_s=\xi$ and $\bar X_s=\bar\xi$.
Moreover, the segment process $\{X(s,\xi)_t, \bar X(s,\bar\xi)_t\}_{t\ge s}$ satisfies
\beq\label{BDD} \bar{\E}^G \sup_{t\in [s,T]} \big(\|X(s,\xi)_t\|_\infty^2+ \|\bar X(s,\bar \xi)_t\|_\infty^2\big)<\infty,\ \ T\in [s,\infty).
\end{equation}
\end{rem}
To characterize the order preservation for solution of \eqref{E1}, we introduce the partial-order on $\C.$ Firstly, for $x=(x^1,\cdots, x^d)$ and $ y=(y^1,\cdots, y^d)\in\R^d$, we write $x\le y$ if $x^i\le y^i$ holds for all $1\le i\le d.$
Similarly,  for $\xi=(\xi^1,\cdots,\xi^d)$ and $\eta=(\eta^1,\cdots,\eta^d)\in\C$, we write $\xi\le \eta$ if $\xi^i(s)\le \eta^i(s)$ holds for all $s\in [-r_0,0]$ and $1\le i\le d.$ %A function $f$ on $\C$ is called increasing if $f(\xi)\le f(\eta)$ for $\xi\le \eta$.
Moreover, for any $\xi_1,\xi_2\in\C$, $\xi_1\land\xi_2\in \C$ is defined by $$(\xi_1\land\xi_2)^i=\min\{\xi^i_1,\xi_2^i\}, \ \  1\le i\le d.$$

\beg{defn} The stochastic differential system  $(\ref{E1})$ is called order-preserving, if  for any $s\ge 0$ and $\xi, \bar\xi\in \C$ with  $\xi\le\bar\xi$,
%$\xi, \bar\xi\in L^2_G(\OO_s\to\C)$ with  $\mathcal{C}$-a.s. $\xi\le\bar\xi$,
it holds  $\mathcal{C}$-q.s.   $$X(s,\xi;t)\le \bar X(s,\bar\xi;t),\ t\ge s.$$   \end{defn}

We first present the following sufficient conditions for the order preservation, which reduce back to the corresponding ones in \cite{HW} when the noise is an $m$-dimensional standard Brownian motion and in \cite{LW} where the system is path independent.

\beg{thm}\label{T1.1} Assume {\bf (H1)}-{\bf (H2)}. The system \eqref{E1} is order-preserving provided that the following  two conditions are satisfied:
\beg{enumerate}
\item[$(1)$] For any  $1\leq i\leq d$, $\xi,\eta\in\C$ with $\xi\leq \eta$ and $\xi^i(0)=\eta^i(0)$,
    \begin{equation*}
b^i(t,\xi)-\bar{b}^i(t, \eta)+2G(h^i(t,\xi)-\bar{h}^i(t, \eta))\leq 0,\ \ \text{a.e.}\ t\ge 0.
    \end{equation*}
\item[$(2)$] For a.e.\ $t\ge 0$ it holds $\si(t, \cdot)= \bar\si(t,\cdot)$    and $\sigma^{ij}(t,\xi)=\si^{ij}(t, \eta)$ for any $1\leq i\leq d$, $1\leq j\leq m$, $\xi,\eta\in \C$ with $\xi^{i}(0)=\eta^i(0)$.
\end{enumerate}
\end{thm}

Condition (2)  means that for a.e. $t\ge 0$, $\si(t,\xi)=\bar\si(t,\xi)$   and  $\si^{ij}(t,\xi)$ only depends on  $t$ and $\xi^i(0)$.

The next result shows that    these conditions are also necessary if all coefficients are continuous on $[0,\infty)\times\C$.

\beg{thm}\label{T1.2} Assume {\bf (H1)}-{\bf (H2)} and that $\eqref{E1}$ is order-preserving. If in addition, $b, h, \si$ and $\bar b, \bar h,\bar \si$ are continuous on $[0,\infty)\times \C$,
then conditions $(1)$ and $(2)$ in Theorem \ref{T1.1} hold.
\end{thm}
These two theorems will be proved in Section 4 and Section 5 respectively.

\section{Proof  of Theorem  \ref{T1.1}}

  Assume  {\bf (H1)}-{\bf (H2)}, and let conditions (1) and (2) hold. For any $T>t_0\ge 0$ and $\xi, \bar\xi\in \C$ with $\xi\le\bar\xi$, it suffices  to prove
  \beq\label{P} \bar{\E}^G\sup_{t\in[t_0,T]} (X^i(t_0,\xi; t)-\bar X^i(t_0,\bar\xi;t))^+=0,\ \ 1\le i\le d,\end{equation} where $s^+:=\max\{0,s\}.$ In fact, by \eqref{P} and \eqref{rep}, for any $P\in\mathcal{P}$, it holds
   \beq\label{P'} \E_P\sup_{t\in[t_0,T]} (X^i(t_0,\xi; t)-\bar X^i(t_0,\bar\xi;t))^+=0,\ \ 1\le i\le d.\end{equation}
   This implies $$P\{X^i(t_0,\xi; t)>\bar X^i(t_0,\bar\xi;t),\ \ t\in[t_0,T]\}=0,$$
   from which we have
   $$\mathcal{C}\{X^i(t_0,\xi; t)>\bar X^i(t_0,\bar\xi;t),\ \ t\in[t_0,T]\}=0.$$
   So the order preservation holds.
  For simplicity, in the following we   denote $X(t)=X(t_0,\xi;t)$ and $\bar X(t)=\bar X(t_0,\bar\xi;t)$ for $t\ge t_0-r_0$. Then it holds
  $$X(t)=\xi(t-t_0),\ \ \bar X(t)= \bar \xi(t-t_0),\ \ t\in [t_0-r_0,t_0].$$

To prove   \eqref{P} using It\^o's formula, we    take  the following $C^2$-approximation of $s^+$  as in the proof of \cite[ Theorem 1.1]{HW}.
For any $n\ge 1$, let $\psi_n: \R\to [0,\infty)$ be constructed as follows: $\psi_n(s)=\psi_n'(s)=0$ for $s\in (-\infty,0]$, and
$$\psi_n''(s)=\beg{cases} 4n^2s, & s\in [0,\ff 1 {2n}],\\
-4n^2(s-\ff 1 n), & s\in [\ff 1 {2n}, \ff 1 n],\\
0, &\text{otherwise}.\end{cases}$$
It is  not difficult to see that
\beq\label{1.3}
0\le \psi_n'\le 1_{(0,\infty)}, \ \text{and\ as\ } n\uparrow\infty: \ 0\le \psi_n(s)\uparrow s^+,\ \ s\psi_n''(s)\le  1_{(0,\ff 1 n)}(s)\downarrow 0.
\end{equation}
%Let $$\tau_k=\inf\big\{t\ge t_0: |X(t)-X(t)\land\bar X(t)|\ge k\big\},\ \ k\ge 1.$$
In view of $$\psi_n(X^i(t_0)-\bar X^i(t_0))=\psi_n(\xi^i(0)-\bar\xi^i(0))=0,$$
and due to (2) $\si(t,\cdot)=\bar\si(t,\cdot)$ for a.e. $t\ge 0$, it follows from It\^o's  formula that
\beq\label{1.4} \beg{split}&\psi_n(X^i(t)-\bar X^i(t))^2\\
&=2\sum_{j=1}^m \int_{t_0}^t (\si^{ij}(s, X_s)- \si^{ij}(s,\bar X_s))\{\psi_n\psi_n'\}(X^i(s)-\bar X^i(s))\d B^j(s)\\
&+2\int_{t_0}^{t} \<h^i(s,X_s)-\bar h^i(s,\bar X_s),\d \<B\>(s)\>\{\psi_n\psi_n'\}(X^i(s)-\bar X^i(s))\d s\\
&+2\int_{t_0}^{t} (b^i(s,X_s)-\bar b^i(s,\bar X_s))\{\psi_n\psi_n'\}(X^i(s)-\bar X^i(s))\d s\\
&\quad +  \sum_{j=1,k=1}^m \int_{t_0}^{t}\{ \psi_n\psi_n''+\psi_n'^2\}(X^i(s)-\bar X^i(s))\\
&\qquad\qquad\qquad\times  (\si^{ij}(s,X_s)- \si^{ij}(s,\bar X_s))(\si^{ik}(s,X_s)- \si^{ik}(s,\bar X_s))\d \<B\>_{jk}(s)\\
 &=M_i(t)+\bar{M}_i(t)+I_1+I_2
 \end{split}\end{equation}for any $n\ge 1$, $1\le i\le d$ and $t\ge t_0,$
where
\begin{align*}
&M_i(t):= 2\sum_{j=1}^m \int_{t_0}^t (\si^{ij}(s, X_s)- \si^{ij}(s,\bar X_s))\{\psi_n\psi_n'\}(X^i(s)-\bar X^i(s))\d B^j(s),\\
&\bar{M}_i(t):=2\int_{t_0}^{t} \<(h^i(s,X_s)-\bar h^i(s,\bar X_s)),\d \<B\>(s)\>\{\psi_n\psi_n'\}(X^i(s)-\bar X^i(s))\d s\\
&\qquad\qquad-4\int_{t_0}^{t} G[\{\psi_n\psi_n'\}(X^i(s)-\bar X^i(s))(h^i(s,X_s)-\bar h^i(s,\bar X_s))]\d s,\\
&I_1:=2\int_{t_0}^{t} (b^i(s,X_s)-\bar b^i(s,\bar X_s))\{\psi_n\psi_n'\}(X^i(s)-\bar X^i(s))\d s\\
&\qquad\qquad+4\int_{t_0}^{t} G[\{\psi_n\psi_n'\}(X^i(s)-\bar X^i(s))(h^i(s,X_s)-\bar h^i(s,\bar X_s))]\d s,\\
&I_2:=\sum_{j=1,k=1}^m \int_{t_0}^{t}\{ \psi_n\psi_n''+\psi_n'^2\}(X^i(s)-\bar X^i(s))\\
&\qquad\qquad\qquad\times  (\si^{ij}(s,X_s)- \si^{ij}(s,\bar X_s))(\si^{ik}(s,X_s)- \si^{ik}(s,\bar X_s))\d \<B\>_{jk}(s).
\end{align*}
Noting that $0\le\psi_n'(X^i(s)-\bar X^i(s))\le 1_{\{X^i(s)>\bar X^i(s)\}}$ and when $X^i(s)>\bar X^i(s)$ one has
$(X_s\land \bar X_s)^i(0)=(\bar X_s)^i(0),$ it follows from (1) that  for a.e. $s\in [t_0,T],$ and $n\geq 1$
\begin{equation}\label{le}
[b^i(s,X_s\land\bar X_s)-\bar b^i(s,\bar X_s)+2G[h^i(s,X_s\land\bar X_s)-\bar h^i(s,\bar X_s)]]\{\psi_n\psi_n'\}(X^i(s)-\bar X^i(s))\leq0.
\end{equation}
In view of $\{\psi_n\psi_n'\}(X^i(s)-\bar X^i(s))\geq 0$, it follows from property (a) of $G$ that \begin{align}\label{Gp}&G[\{\psi_n\psi_n'\}(X^i(s)-\bar X^i(s))(h^i(s,X_s)-\bar h^i(s,\bar X_s))]\\ \nonumber
&=\{\psi_n\psi_n'\}(X^i(s)-\bar X^i(s))G(h^i(s,X_s)-\bar h^i(s,\bar X_s)).
\end{align}
For simplicity, let $\Phi^n_s=\{\psi_n\psi_n'\}(X^i(s)-\bar X^i(s))$. Combining \eqref{le} with \eqref{Gp}, {\bf (H1)}, $0\le\psi_n'\le 1$ and properties (b) and (c) of $G$, we obtain
\beg{align*}
&I_1=2\int_{t_0}^{t}[b^i(s,X_s)-\bar b^i(s,\bar X_s)+2G(h^i(s,X_s)-\bar h^i(s,\bar X_s))]\Phi^n_s\d s\\
& \leq 2\int_{t_0}^{t}[b^i(s,X_s)-b^i(s,X_s\land \bar X_s)+2G(h^i(s,X_s)-h^i(s,X_s\land \bar X_s))]\Phi^n_s\d s\\
&+2\int_{t_0}^{t} [b^i(s,X_s\land\bar X_s)-\bar b^i(s,\bar X_s)+2G(h^i(s,X_s\land\bar X_s)-\bar h^i(s,\bar X_s))]\Phi^n_s\d s\\
&\leq2\int_{t_0}^{t}[b^i(s,X_s)-b^i(s,X_s\land \bar X_s)+2G(h^i(s,X_s)-h^i(s,X_s\land \bar X_s))]\Phi^n_s\d s\\
&\leq\int_{t_0}^{t}[b^i(s,X_s)-b^i(s,X_s\land \bar X_s)+2G(h^i(s,X_s)-h^i(s,X_s\land \bar X_s))]^2\d s\\
&\qquad +\int_{t_0}^{t}\psi_n(X^i(s)-\bar X^i(s))^2\d s\\
&\leq\int_{t_0}^{t}C(T,\bar{\sigma})\|X_s-X_s\land \bar X_s\|_{\infty}^2\d s+\int_{t_0}^{t}\psi_n(X^i(s)-\bar X^i(s))^2\d s,\ \ n\geq1,t\in[t_0,T].
\end{align*}
Next, by condition (2) in Theorem \ref{T1.1}, for a.e. $s\in [t_0,T]$, $\si^{ij}(s, X_s)=\bar\si^{ij}(s, X_s)$ and $\si^{ij}(s, X_s)$ depends only on $s$ and $X^i(s)$. So,   \eqref{1.3}, {\bf (H1)} and the positive definite property of $\<B\>$ yield
\beq\label{1.6} \beg{split}
&I_2\leq  \int_{t_0}^{t} \left(1_{\{X^i(s)-\bar X^i(s)\in (0,\ff 1 n)\}}+1_{\{X^i(s)-\bar X^i(s)\in (0,\infty)\}}\right)\\
 &\qquad\qquad \times \sum_{j=1,k=1}^m(\si^{ij}(s,X_s)- \si^{ij}(s,\bar X_s))(\si^{ik}(s,X_s)- \si^{ik}(s,\bar X_s))\d \<B\>_{jk}(s) \\
&\leq  \int_{t_0}^{t} C(T,\bar{\sigma})\{(X^i(s)- \bar X^i(s))^+\}^2\d s,\ \ n\ge 1, t\in [t_0,T].
\end{split}\end{equation}
By the Burkholder-Davis-Gundy  inequality in \cite[Theorem 2.1]{G} or \cite[Lemma 8.1.12]{peng4}, we deduce
\beg{align*}  &\bar{\E}^G\sup_{s\in [t_0,t]}  M_i(s)\\
&\leq C(T,\bar{\sigma})\bar{\E}^G\bigg\{\int_{t_0}^{t}   \big|(\psi_n\psi_n')(X^i(s)-\bar X^i(s))\big|^2 \sum_{j=1}^m|\si^{ij}(s,X_s)- \si^{ij}(s,\bar X_s)|^2\d s\bigg\}^{\frac{1}{2}}\\
 & \le C(T,\bar{\sigma}) \bar{\E}^G \bigg(\int_{t_0}^{t}  \{(X^i(s)- \bar X^i(s))^+\}^2  \psi_n(X^i(s)-\bar X^i(s))^2\d s\bigg)^{\ff 1 2}\\
& \le C(T,\bar{\sigma}) \bar{\E}^G\int_{t_0}^{t} \|X_s-X_s\land \bar X_s\|_\infty^2 \d s\\
 &\quad+  \ff 1 8 \bar{\E}^G\sup_{s\in [t_0, t]}\psi_n(X^i(s)-\bar X^i(s))^2,\ \ n\ge 1, t\in [t_0,T].
\end{align*}
Finally, since $\bar{M}_i$ is a non-increasing $G$-martingale, we obtain from \eqref{Gp} that  $$\bar{\E}^G\sup_{s\in [t_0,t]}\bar{M}_i(s)\leq  0.$$
Now, letting
$$\phi(s)=\sup_{r\in [t_0-r_0, s]}|X(r)-X(r)\land \bar X(r)|^2,\ \ s\in [t_0,T],$$
we obtain
\beq\label{1.8} \beg{split}
&\bar{\E}^{G}\sup_{r\in [t_0-r_0, t]} \psi_n(X^i(r)-\bar X^i(r))^2=\bar{\E}^G\sup_{r\in [t_0, t]} \psi_n(X^i(r)-\bar X^i(r))^2\\
&\leq C\bar{\E}^{G}\int_{t_0}^{t}\|X_{s}-X_{s}\land \bar X_{s}\|_{\infty}^2\d s+\ff 1 8 \bar{\E}^G\sup_{s\in [t_0, t]}\psi_n(X^i(s)-\bar X^i(s))^2,
\end{split}\end{equation}
for some constant $C>0$ and all $n\ge 1, t\in [t_0,T], 1\le i\le d.$ Therefore, for any $n\ge 1$ and $t\in [t_0,T],$ it holds
\beg{equation*}\beg{split} &  \sum_{i=1}^d\bar{\E}^G\sup_{r\in [t_0-r_0, t]} \psi_n(X^i(r)-\bar X^i(r))^2\le C\int_{t_0}^{t}\bar{\E}^G\phi(s)\d s,\ \ n\ge 1.\end{split}\end{equation*}
 Letting  $n\uparrow\infty$, by the monotone convergence theorem in \cite[Theorem 6.1.14]{peng4}, we arrive at
$$\bar{\E}^G  \phi(t)\leq \sum_{i=1}^d\bar{\E}^G\sup_{r\in [t_0-r_0, t]}\{(X^i(r)-\bar X^i(r))^+\}^2\le C  \int_{t_0}^t  \bar{\E}^G\phi(s)\d s,\ \ t\in [t_0,T].$$ By the definition of $\phi$ and \eqref{BDD}, Gronwall's inequality  implies
$$ \bar{\E }^G\phi(T)=0.$$ Thus, we  prove  (\ref{P}).

\section{Proof of Theorem \ref{T1.2}}
\begin{proof}[$\textbf{Proof of (1)}$] Let $1\leq i\leq d$ be fixed.
For any $t_0\ge 0$ and $\xi, \eta\in \C$ with $\xi\le\eta$ and $\xi^i(0)=\eta^i(0)$, it holds $\mathcal{C}$-q.s.
\begin{align}\label{com}X(t_0,\xi;t)\leq \bar{X}(t_0,\eta;t),\ \ t\geq t_0.\end{align}
For simplicity, let $X(t)=X(t_0,\xi;t)$ and $\bar X(t)=\bar X(t_0,\eta;t)$ for $t\ge t_0-r_0$.
For any $\gamma\in \mathbb{S}_+^m\bigcap[\underline{\sigma}^2\textbf{I}_{m\times m}, \bar{\sigma}^2\textbf{I}_{m\times m}]$, take $\theta_s=\sqrt{\gamma},s\geq 0$ and denote $\E_{\P_\theta}=\E_{\P_\gamma}$. Then $\P_{\gamma}$-a.s. $\<B\>(r)=r\gamma$.
By \eqref{E1}, \eqref{rep2} and \eqref{com}, for any $s\ge 0$, we obtain $\P_\gamma$-a.s.
\beq\label{XA}\beg{split}
0&\ge X^i(t_0+s)-\bar X^i(t_0+s)   =\xi^i(0)-\eta^i(0)\\
 &\quad +   \int_{t_0}^{t_0+s}[b^i(r,X_r)-\bar b^i(r,\bar X_r)]\,\d r+\int_{t_0}^{t_0+s}\<h^i(r,X_r)-\bar h^i(r,\bar X_r),\d\<B\>(r)\> \\
&\quad +  \sum_{j=1}^m \int_{t_0}^{t_0+s}[\sigma^{ij}(r,X_r)-\bar \sigma^{ij}(r,\bar X_r)]\,\d B^j(r)\\
&=\int_{t_0}^{t_0+s}[b^i(r,X_r)-\bar b^i(r,\bar X_r)]\,\d r+\int_{t_0}^{t_0+s}\<h^i(r,X_r)-\bar h^i(r,\bar X_r),\gamma\>\d r \\
&\quad +  \sum_{j=1}^m \int_{t_0}^{t_0+s}[\sigma^{ij}(r,X_r)-\bar \sigma^{ij}(r,\bar X_r)]\,\d B^j(r).
\end{split}\end{equation}
By {\bf(H1)}, {\bf(H2)} and \eqref{BDD}, taking expectation in \eqref{XA} under $\P_\gamma$,   we
obtain
\begin{equation}\begin{split}\label{sf}
 \frac{1}{s}\int_{t_0}^{t_0+s}\E_{\P_\gamma}\{[b^i(r,X_r)-\bar b^i(r,\bar X_r)]+\< h^i(r, X_r)-\bar{h}^i(r,\bar{X}_r),\gamma\>\}\d r\leq 0,\ \ s>0.
\end{split}\end{equation}
Thus, taking $s\downarrow 0$ in \eqref{sf}, it follows from \eqref{BDD}, \eqref{rep2}, the continuity of $b,\bar{b},h,\bar{h}$ and dominated convergence theorem that
\begin{equation*}\begin{split}
[b^i(t_0,\xi)-\bar b^i(t_0,\eta)]+\<h^i(t_0,\xi)-\bar h^i(t_0,\eta),\gamma\>\leq 0.
\end{split}\end{equation*}
By the definition of $G$, we derive
\begin{equation*}\begin{split}
&[b^i(t_0,\xi)-\bar b^i(t_0,\eta)]+2G(h^i(t_0,\xi)-\bar h^i(t_0,\eta))\\
&=[b^i(t_0,\xi)-\bar b^i(t_0,\eta)]+2\sup_{\gamma\in \mathbb{S}_+^m\bigcap[\underline{\sigma}^2\textbf{I}_{m\times m}, \bar{\sigma}^2\textbf{I}_{m\times m}]}\frac{1}{2}\<h^i(t_0,\xi)-\bar h^i(t_0,\eta),\gamma\>\leq 0.
\end{split}\end{equation*}
The proof is completed.
\end{proof}

\beg{proof}[\textbf{Proof of (2)}] For any $t_0\geq 0$ and $\xi,\bar{\xi}\in\C$ with $\xi\leq \bar{\xi}$, it holds $\mathcal{C}$-q.s.
$$X(t_0,\xi;t)\leq \bar{X}(t_0,\bar{\xi};t), \ \ t\geq t_0.$$
Taking $\theta_s=\bar{\sigma}$, \eqref{rep2} implies $\P_\theta$-a.s.
$X(t_0,\xi;t)\leq \bar{X}(t_0,\bar{\xi};t), \ \ t\geq t_0.$ Noting that $\P_\theta$-a.s. $\<B\>(r)=\bar{\sigma}^2r$, \eqref{E1} reduces to the SDE driven by Brownian motion under $\P_\theta$. According to the necessary condition of order preservation for functional SDEs in \cite[Theorem 1.2 (II)]{HW}, we immediately get the results desired.
\end{proof}

\end{document}